
\magnification\magstep 1
\baselineskip=0,59 true cm
\vsize=21 true cm
\topinsert\vskip 2 true cm
\endinsert
{\centerline{\bf{LINEAR CONNECTIONS ON THE TWO PARAMETER QUANTUM PLANE}}}
\vskip 2 true cm
{\centerline{\bf{Yvon Georgelin$^a$, Thierry Masson$^b$, Jean-Christophe
Wallet$^a$}}}
\vskip 1 true cm
{\centerline{$^a$ Division de Physique Th\'eorique{\footnote
\dag{\sevenrm{Unit\'e de Recherche des Universit\'es Paris 11 et Paris 6
associ\'ee au CNRS}}}, Institut de Physique
Nucl\'eaire}}
{\centerline{F-91406 ORSAY Cedex, France}}
\vskip 0,5 true cm
{\centerline{$^b$ Laboratoire de Physique Th\'eorique et Hautes
Energies{\footnote*{\sevenrm{Laboratoire associ\'e au CNRS, URA D0063}}} }}
{\centerline{Universit\'e de Paris-Sud, B\^at 211, F-91405 ORSAY, France}}
\vskip 1,5 true cm
{\bf{Abstract:}} We apply a recently proposed definition of a linear connection
in non commutative geometry based on the natural bimodule structure of the
algebra of differential forms to the case of the two-parameter quantum plane.
We find that there exists a non trivial family of linear connections only when
the two parameters obeys a specific relation.\par
\vskip 3 true cm
{\noindent{IPN}}O-TH-9531 (May 1995)\hfill\break
LPTHE Orsay 95/43
\vfill\eject
{\bf{I) Introduction}}\hfill\break
\vskip 0,5 true cm
In the last few years, much attention has been attracted by non commutative
differential geometry. Many related attempts focused on generalization of
differential forms [1,2,3], as well as covariant derivative [1].
This later generalization of covariant derivative only used a left or a right
module structure. However, in order to extend the notion of linear connection
to non commutative geometry, one has to deal with one-forms. Then, the bimodule
structure of the space of one-forms should be taken into account. This
has been done in
[4a,b] for the general derivation based differential calculus
and in [5] for more general differential calculi. Other examples based on [5]
have been worked out in [6] and [7]. \par
In this letter, we adopt this last viewpoint to construct linear connections on
the two parameter quantum plane. In section II, we recall the general
definition of a linear connection on a non commutative algebra as well as some
results already obtained in [6] on the one parameter quantum plane
[8]. In section III we present our results for the two parameter quantum
plane.\par
\vskip 0,2 true cm
{\bf{II) Linear connections over a noncommutative algebra}}\par
{\bf{IIa) Basic tools}}\par
\vskip 0,5 true cm
Let ${\cal{A}}$ and ($\Omega^*({\cal{A}}),d$) be respectively a noncommutative
algebra and a differential calculus over it ($d$ is the exterior
derivative). Let $\Omega^k({\cal{A}})$ ($k\ge 0$) be the algebra of
differential forms of degree $k$ on ${\cal{A}}$ and
$\pi\colon\ \Omega^1({\cal{A}})\otimes_{\cal{A}}\Omega^1({\cal{A}})\to
\Omega^2({\cal{A}})$ the projection defined
by the product of forms ($\otimes_{\cal{A}}$ is the tensor product on
${\cal{A}}$).\par
We first recall the definition of a linear connection on ${\cal{A}}$ given in
[4-6] that will be the main ingredient used in this letter.\par
{\bf{Definition:}} A linear connection over ${\cal{A}}$ is determined from two
maps $\sigma$ and $D$, \hfill\break
$\sigma\colon\ \Omega^1({\cal{A}})\otimes_{\cal{A}}\Omega^1({\cal{A}})\to
\Omega^1({\cal{A}})\otimes_{\cal{A}}\Omega^1({\cal{A}})$,
$D\colon\ \Omega^1({\cal{A}})\to
\Omega^1({\cal{A}})\otimes_{\cal{A}}\Omega^1({\cal{A}})$,
satisfying the following properties: for any
$f,g\in\Omega^0({\cal{A}})$, $\alpha,\beta\in\Omega^1({\cal{A}})$\par
$${\hbox{$\sigma$ is left and right--linear}}\
\sigma(f\alpha\otimes\beta g)=f\sigma(\alpha\otimes\beta) g
\eqno(2.1)$$
$$\pi(\sigma+1)=0 \eqno(2.2)$$
$$D(f\alpha)=df\otimes\alpha+fD\alpha,\
D(\alpha f)=D\alpha f+\sigma(\alpha\otimes df)\eqno(2.3).$$
Concerning this definition, some comments are in order. \par
Basically, the map $\sigma$, which is a
bimodule automorphism in $\Omega^1({\cal{A}})\otimes_{\cal{A}}
\Omega^1({\cal{A}})$,
can be viewed as a generalization of the
permutation map that would appear in the case of a commutative algebra. The
map $D$ is a
noncommutative extension of the covariant derivative (observe in particular the
Leibnitz rules, property (2.3)). Furthermore, it can be easily
seen that the right (resp. left)
${\cal{A}}$--linearity of $\sigma$, property (2.1), insures that $D((\xi f)g)=
D(\xi(fg))$ (resp. $D(d((fg)h))=D(d(f(gh)))$ ) for any
$\xi\in\Omega^1({\cal{A}})$, $f,g,h\in\Omega^0({\cal{A}})$.\par
As far as the properties (2.2) and (2.3) are
concerned, it must be pointed out that
not each solution of (2.2) admits a covariant derivative, that is, a map $D$
fullfilling the Leibnitz rules (2.3). In other words, a linear
connection for each
$\sigma$ does not necessarily exist.\par
Finally, the above definition is nothing but a noncommutative extension of the
definition of a linear connection in term of a covariant derivative introduced
by Koszul [9] in the framework of commutative geometry. Indeed, for a
commutative algebra ${\cal{A}}$, it follows from (2.3) that $\sigma$
reduces to the usual permutation, so that
 (2.1) and (2.2) are satisfied. The Koszul definition then follows.\par
\vskip 0,2 true cm
{\bf{IIb) Application to the Manin quantum plane}}\par
\vskip 0,5 true cm
Let $\Omega^*=\Omega^0\oplus\Omega^1\oplus\Omega^2$ the algebra
of forms of the Manin
quantum plane [8] whose generators $x^i\colon=(x,y)$,
$\xi^i=dx^i\colon=(\xi,\eta)$ obey the following relations
$$x^ix^j-q^{-1}\hat R^{ij}_{kl}x^kx^l=0;\ x^i\xi^j-
q\hat R^{ij}_{kl}\xi^kx^l=0;\ \xi^i\xi^j+
q\hat R^{ij}_{kl}\xi^k\xi^l=0  \eqno(2.4a;b;c),$$
where the tensor $\hat R^{ij}_{kl}$ is given in matrix form by
$$\hat R=\pmatrix{q&0&0&0\cr
0&q-q^{-1}&1&0\cr
0&1&0&0\cr
0&0&0&q\cr}   \eqno(2.5).$$
We recall that $\Omega^*$ is stable under the action of $SL_q(2,C)$ since the
following
identity $\hat R^{ij}_{kl}a^k_ma^l_n=$$a^i_ka^j_l\hat R^{kl}_{mn}$ holds (where
$a^i_j$ denote generically the four generators of $SL_q(2,C)$ written again in
matrix form with $da^i_j=0$ and $d$ is the exterior derivative).\par
Now, the action of $D$ on eqn. (2.4b) which is defined on
$\Omega^1$ determines the action of $\sigma$ on
$\xi^i\otimes\xi^j$. It is given [6] by
$$\sigma(\xi\otimes\xi)=q^{-2}(\xi\otimes\xi);\
\sigma(\xi\otimes\eta)=q^{-1}(\eta\otimes\xi) \eqno(2.6a;b),$$
$$\sigma(\eta\otimes\xi)=q^{-1}(\xi\otimes\eta)+(q^{-2}-1)(\eta\otimes\xi);\
\sigma(\eta\otimes\eta)=q^{-2}(\eta\otimes\eta)  \eqno(2.6c;d),$$
so that $\sigma=q^{-1}\hat R^{-1}$ and can be proven to be left and right
linear. Furthermore, $\sigma$ obeys the property (2.2) of the general
definition
, namely $\pi(\sigma+1)=0$, as it can be easily seen by computing the action of
$\pi(\sigma+1)$ on $\xi^i\otimes\xi^j$ and using the relations (2.4c). Notice
that $\sigma^2\ne1$; however, it satisfies a Hecke relation given by
$$(\sigma+1)(\sigma-q^{-2})=0  \eqno(2.7),$$
where the simple (resp. triply degenerate) eigenvalue $-1$ (resp. $q^{-2}$)
corresponds to the antisymmetric (resp. symmetric) eigenspaces with
eigenvectors $\xi\otimes\eta-q\eta\otimes\xi$ (resp. $\xi\otimes\xi$,
$\eta\otimes\eta$, $\eta\otimes\xi+q\xi\otimes\eta$).\par
Finally, the covariant derivative map $D$ acting on the $\xi^i$'s can be cast
into the form
$$D\xi^i=\mu\ x^i\theta\otimes\theta  \eqno(2.8),$$
where $\mu$ is an arbitrary parameter (which has the dimension of a mass) and
$\theta\in\Omega^1$ is the unique 1-form invariant under the coaction of
$SL_q(2,C)$ which is given by
$\theta=x\eta-qy\xi$ (up to an overall constant) with $\theta^2=0$.\par
\vskip 0,2 true cm
{\bf{III) The two parameter quantum plane}}\par
\vskip 0,5 true cm
The algebraic structure of the two parameter quantum plane is now defined by
$$x^ix^j-q^{-1}\hat R^{ij}_{kl}(p,q)x^kx^l=0;\ x^i\xi^j-
p\hat R^{ij}_{kl}(p,q)\xi^kx^l=0;\ \xi^i\xi^j+
p\hat R^{ij}_{kl}(p,q)\xi^k\xi^l=0  \eqno(3.1a;b;c),$$
where again $x^i\colon=(x,y)$, $\xi^i=dx^i\colon=(\xi,\eta)$ and the tensor
$\hat R^{ij}_{kl}(p,q)$ is given by
$$\hat R=\pmatrix{q&0&0&0\cr
0&q-p^{-1}&qp^{-1}&0\cr
0&1&0&0\cr
0&0&0&q\cr}   \eqno(3.2),$$
and reduces to (2.5) when $q=p$. Recall that $\hat R(p,q)$ satisfies the
braid equation: \hfill\break
$\hat R_{12}(p,q)\hat R_{23}(p,q)\hat R_{12}(p,q)$$=
\hat R_{23}(p,q)\hat R_{12}(p,q)\hat R_{23}(p,q)$ with as usual $
\hat R_{12}(p,q)\colon=\hat R(p,q)\otimes 1$, $
\hat R_{23}(p,q)\colon=1\otimes\hat R(p,q)$. For the
moment, $p$ and $q$ are independant arbitrary
parameters. We will see in a while that a non trivial family of linear
connections can be obtained only when $p$ and $q$ obeys a supplementary
relation.\par
{}From the action of $D$ on eqn.(3.1b), we determine the action of $\sigma$ on
$\xi^i\otimes\xi^j$. Namely, using the property (2.3) of the general
definition, we obtain generically
$$\xi^i\otimes\xi^j=p\hat R^{ij}_{kl}(p,q)\sigma(\xi^k\otimes\xi^l);\
x^iD\xi^j=p\hat R^{ij}_{kl}(p,q)(D\xi^k)x^l  \eqno(3.3a;b).$$
The equation (3.3a) is verified when $\sigma=p^{-1}\hat R^{-1}(p,q)$ or
equivalently
$$\sigma(\xi\otimes\xi)=p^{-1}q^{-1}(\xi\otimes\xi);\
\sigma(\xi\otimes\eta)=p^{-1}(\eta\otimes\xi)  \eqno(3.4a;b),$$
$$\sigma(\eta\otimes\xi)=q^{-1}(\xi\otimes\eta)
+(p^{-1}q^{-1}-1)(\eta\otimes\xi);\ \sigma(\eta\otimes\eta)=p^{-1}q^{-1}
(\eta\otimes\eta)  \eqno(3.4c;d).$$
It is easy to verify the left and right-linearity of $\sigma$ (property (2.1)).
Notice that right-linearity stemms from left-linearity, as a mere consequence
of the braid equation for $\hat R(p,q)$. To see that, it is
sufficient to combine $\sigma=p^{-1}\hat R^{-1}(p,q)$ with (3.3a) and (3.1b)
together with the braid equation; then, the above
statement follows. Notice also that $\sigma$ fullfills a braid equation
since $\hat R(p,q)$ does. Besides, combining (3.4) and (3.1c), we also
verify that
$\pi(\sigma+1)=0$ (property (2.2)) still holds.\par
Thus, we have determined a suitable $\sigma$ map on the two parameter quantum
plane. The remaining eqn. (3.3b) will fix a relation between $p$ and $q$
so that a non trivial family of covariant derivative $D$ can be associated to
this $\sigma$. By non trivial family, we mean that $D$ can
differ from the trivial case $D\xi=D\eta=0$ which always verifies
(3.3b). Namely, we find after
some calculations that (3.3b) is verified provided $D\xi$ and $D\eta$ are given
by
$$D\xi=x Z\theta\otimes\theta;\ D\eta=yZ
\theta\otimes\theta;\ Z=\mu x^{n-1}y^{n-1}  \eqno(3.5a;b;c),$$
where $\mu$ is an overall complex parameter; (3.5) must be supplemented by
$$p=q^n\ \ n\ge 1  \eqno(3.6).$$
In (3.5), $\theta=x\eta-qy\xi$ and still verifies $\theta^2=0$ and
$$\sigma(\xi\otimes\theta)=q^{-1-2n}\theta\otimes\xi;\
\sigma(\theta\otimes\xi)=q^n(\xi\otimes\theta)-
(1-q^{-1-n})(\theta\otimes\xi)
  \eqno(3.7a;b),$$
$$\sigma(\eta\otimes\theta)=q^{-2-n}\theta\otimes\eta;\
\sigma(\theta\otimes\eta)=q(\eta\otimes\theta)-(1-q^{-1-n})
(\theta\otimes\eta)  \eqno(3.7c;d),$$
$$\sigma(\theta\otimes\theta)=q^{-1-n}(\theta\otimes\theta)  \eqno(3.7e),$$
where we used (3.6).\par
It is interesting to observe that a non trivial family of
linear connection on the two parameter
quantum plane can be consistently found only when the two parameters are
related each other through (3.6).\par
Therefore, the corresponding map $\sigma$ is given in
the tensor form by
$$\sigma=\pmatrix{q^{-1-n}&0&0&0\cr
0&0&q^{-n}&0\cr
0&q^{-1}&q^{-1-n}-1&0\cr
0&0&0&q^{-1-n}\cr}  \eqno(3.8)$$
and satisfies the Hecke relation $(\sigma+1)(\sigma-q^{-1-n})=0$
where the simple (resp. triply degenerate) eigenvalue $-1$ (resp. $q^{-1-n}$)
corresponds to the antisymmetric (resp. symmetric) eigenspace with
eigenvectors $\xi\otimes\eta-q^n\eta\otimes\xi$ (resp. $\xi\otimes\xi$,
$\eta\otimes\eta$, $\eta\otimes\xi+q\xi\otimes\eta$).\par
Some remarks are in order. Firstly, it happens that $\sigma$ is actually the
only generalized permutation for which there exits a covariant
derivative.\hfill\break
Next, one recovers the results of section IIb
when $n=1$. We point out that the case $p=q$ is very similar, as far as
the structure of the set of linear connections is concerned, to the cases
$p=q^n$. \hfill\break
Finally, in (3.6) we have restricted $n\ge 1$ since working on the
quantum plane forces to consider only positive powers in $x$ and $y$. However,
the case $pq=1$ is interesting. In this last situation, the algebra is
formally the non commutative torus (for which negative powers of $x$ and $y$
are allowed). The differential
calculus is then based on derivations in the sense of [4b]. An easy
calculation using $\sigma$ given in (3.4) for $pq=1$
leads to an eight complex parameter family of linear connections. This family
can also be obtained using the definition of linear connection based on
derivations given in [2,4a].\par
Let $\pi_{12}\colon=\pi\otimes1$ and $\sigma_{12}\colon=\sigma\otimes1$. The
covariant derivative map defined above can be extended to
$\Omega^1\otimes_{\cal{A}}\Omega^1$ by
$$D(\alpha\otimes\beta)=D\alpha\otimes\beta+\sigma_{12}(\alpha\otimes D\beta)
\eqno(3.9),$$
for any $\alpha,\beta\in\Omega^1$. Consider now the following map
$$\pi_{12}D^2\colon\Omega^1\to(\Omega^2/\Theta)\otimes_{\cal{A}}\Omega^1
\eqno(3.10),$$
where $\Theta$ is a submodule of $\Omega^2$, called the torsion module, given
by the image of $\Omega^1$ by $d-\pi D$. This map is left-linear, namely
$\pi_{12}D^2(f\alpha)=f\pi_{12}D^2(\alpha)$ for any $f\in\Omega^0$ ($\Omega^0=
{\cal{A}}$),
$\alpha\in\Omega^1$ since $\sigma$ verifies $\pi(\sigma+1)=0$. \par
By noticing that $(d-\pi D)\xi^i=0$ holds, thanks to $\theta^2=0$ and combining
(3.1), (3.9) and $\pi(\sigma+1)=0$, we obtain
$$\pi_{12}D^2\theta=0  \eqno(3.11).$$
Now, using again (3.1), (3.5), (3.6), (3.8) we find that
$$\pi_{12}D^2\xi^i=\Omega^i\otimes\theta  \eqno(3.12),$$
with
$$\Omega^i=f(q)x^iZ\xi\eta  \eqno(3.13a),$$
$$f(q)={{1}\over{1-q^{n+1}}}\big(1-q^{3(n+1)}-q^{n+1}+q^{(2-n)(n+1)}-
q^{(1-n^2)}+q^{(n+1)(1-2n)}\big) \eqno(3.13b).$$
Then, from (3.12) and (3.13), we obtain the 2-form curvature $\Omega^i_j$
given by
$$\pi_{12}D^2\xi^i=-\Omega^i_j\otimes\xi^j  \eqno(3.14),$$
$$\Omega^i_j=f(q)\pmatrix{q^{-n-1}xy&-q^{-1-2n}x^2\cr
q^{-1-n}y^2&-q^{-1-2n}yx\cr}Z\xi\eta   \eqno(3.15).$$
\vskip 0,2 true cm
{\bf{IV) Conclusion}}\par
\vskip 0,5 true cm
In this letter, working on the two parameter quantum plane, we have
shown that there exists a non trivial family of linear connections only when
the parameters are related through $p=q^n$, $n\ge 1$. In this respect, the
usual one parameter quantum plane ($n=1$) is fully representative of
the $p=q^n$ situation. It remains to see whether this specific
relation between $p$ and $q$ is purely accidental or reflects a deeper
property.
\vskip 2 true cm
{\bf{Acknowledgement:}}\par
We are very grateful to M. Dubois-Violette and J. Madore for numerous
discussions and comments.\par
\vfill\eject
{\bf{REFERENCES}}\par
\vskip 2 true cm
\item{[1]:} Connes A., Non-commutative differential geometry, Inst. Hautes
Etudes Sci. Publ. Math. 62 (1986) 257; see also Connes A., G\'eom\'etrie non
commutative, InterEditions, Paris, 1990.
\item{[2]:} Dubois-Violette M., C. R. Acad. Sci. Paris 307, S\'erie I (1988)
403.
\item{[3]:} Wess J., Zumino B., Nucl. Phys. B(Proc. Suppl.)18 (1990) 302.
\item{[4a]:} Dubois-Violette M., Michor P., Connections in Central Bimodules,
Orsay LPTHE preprint 94/100 (1994).
\item{[4b]:} Dubois-Violette M., Michor P., C. R. Acad. Sci. Paris 319, S\'erie
I (1994) 927.
\item{[5]:} Mourad J., Class. Quant. Grav. 21 (1995) 965.
\item{[6]:} Dubois-Violette M., Madore J., Masson T., Mourad J., Orsay LPTHE
preprint 94/94 to appear in Lett. Math. Phys.
\item{[7]:} Madore J., Masson T., Mourad J., Orsay LPTHE preprint 94/96 (1994)
to appear in Class. Quant. Grav.
\item{[8]:} Manin Yu., Quantum groups and non comuutative geometry, in Pub. CRM
1561, University of Montreal, 1988.
\item{[9]:} Koszul J.L., Lectures on Fibre Bundles and Differential Geometry,
Tata Institute of Fundamental Research, Bombay (1960).

\end